\documentclass[11pt]{amsart}

\usepackage{amssymb,amsmath,amsthm,amstext,amsfonts}
\usepackage{epic,eepic}
\usepackage[dvips]{graphicx}
\usepackage{psfrag}
\usepackage{txfonts}
\usepackage{color}
\usepackage{mathrsfs}
\usepackage{upgreek}
\usepackage{epstopdf}
\usepackage{caption}
\usepackage{subcaption}
\usepackage{bbm}

\makeatletter \@addtoreset{equation}{section} \makeatother

\theoremstyle{plain}

\newtheorem{theorem}{Theorem }[section]
\newtheorem{proposition}[theorem]{Proposition}
\newtheorem{lemma}[theorem]{Lemma}

\theoremstyle{definition} \theoremstyle{remark}
\newtheorem{remark}[theorem]{Remark}

\DeclareMathAlphabet{\mathpzc}{OT1}{pzc}{m}{it}

\newcommand{\Fix}{\textnormal{Fix}}

\newcommand{\SO}{\textnormal{SO}}
\newcommand{\Z}{\mathbb{Z}}
\newcommand{\R}{{\mathbb R}}
\newcommand{\D}{{\mathbb D}}

\begin{document}

\title{A hybrid heteroclinic cycle}

\author{Sofia B.S.D.\ Castro}
\address[S.\ Castro]{Faculdade de Economia and Centro de Matem\'atica, Universidade do Porto, Rua Dr.\ Roberto Frias, 4200-464 Porto, Portugal.}
\email{sdcastro@fep.up.pt}

\author{Alexander Lohse}
\address[A.\ Lohse]{Department of Mathematics, Universit\"at Hamburg, Bundesstra{\ss}e 55, 20146 Hamburg, Germany.}
\email[Corresponding author]{alexander.lohse@math.uni-hamburg.de}

\begin{abstract}
Using a vector field in $\R^4$, we provide an example of a robust heteroclinic cycle between two equilibria that displays a mix of features exhibited by well-known types of low-dimensional heteroclinic structures, including simple, quasi-simple and pseudo-simple cycles. Our cycle consists of two equilibria on one coordinate axis and two connections. One of the connections is one-dimensional while the other is two-dimensional.
We compare our heteroclinic cycle to others in the literature that are similar in architecture, and illustrate how the standard methods used to analyse those cycles fail to provide sufficient information on the attraction properties of our example.
 The instability of two subcycles contained in invariant three-dimensional subspaces seems to indicate that our cycle is generically completely unstable.
 Although this cycle is one of the simplest possible and exists in low-dimension, the complete study of the stability of our cycle by using the standard techniques for return map reduction is not possible given the hybrid nature of the return map. 
\end{abstract}

\maketitle


\section{Introduction}

A heteroclinic connection is a solution to a dynamical system $\dot x = f(x)$ that lies in the intersection of the stable and unstable manifold of two invariant sets, usually equilibria. Heteroclinic dynamics appear in various real-life systems with intermittent behaviour, ranging from fluid dynamics to Lotka-Volterra-type models. Their robustness with respect to perturbations of the system usually comes as a consequence of the presence of invariant subspaces, which is often -- though not always -- induced by a group action. This makes them interesting in their own right, as potential robust attractors in equivariant dynamical systems.
The symmetry group constrains the possible heteroclinic connection structures in a system allowing for a sharp distinction among heteroclinic cycles with the same architecture, but arising in systems with different symmetry groups.
This becomes evident e.g.\ in the study of so-called simple heteroclinic cycles in $\R^4$, which are classified into types A, B, C and Z, depending on properties of the symmetry group, see \cite{CastroLohse2016, PodviginaLohse}.

The stability of heteroclinic cycles is of interest since it determines whether or not the dynamics they prescribe are observable in experiments and simulations. When a heteroclinic cycle can be classified as either simple, quasi-simple or pseudo-simple cycles, there are various results that enable a systematic study of stability. See \cite{KrupaMelbourne2004}, \cite{PodviginaAshwin2011} for the stability of simple cycles, \cite{GarridoDaSilvaCastro} for that of quasi-simple cycles and \cite{PodviginaChossat} for pseudo-simple cycles.

One of the most elementary heteroclinic objects is a cycle with two hyperbolic equilibria and connections in both directions between these equilibria. Its architecture is that of an oriented graph with two nodes and one edge in each direction. It appears in equivariant systems with different symmetry groups, displaying a variety of stability configurations, see \cite{CastroLohse2014, CastroLohse2016, ChossatLohsePodvigina}. In this paper we construct a class of equivariant vector fields in $\R^4$ supporting such a heteroclinic cycle. There is a one-dimensional connection from one equilibrium to the other and a two-dimensional connection in the opposite direction. Our construction is a modification of an example in Castro and Lohse \cite{CastroLohse2014}. We show that the global transitions in our example possess different properties: one is reminiscent of a type A transition, the other of a type B transition. This makes it impossible to use standard methods, such as those used in \cite{KrupaMelbourne2004}, \cite{PodviginaAshwin2011} or \cite{GarridoDaSilvaCastro}, to draw conclusions about the stability in this system. This observation emphasizes the variety of heteroclinic cycles whose connection structure can be described by the same, very simple, graph. It also exposes the complexity attached to the study of stability of heteroclinic objects.

We finish this section with a brief description of the essential definitions and concepts. The following section constructs the class of vector fields which support a heteroclinic cycle that does not fit into the categories above. The final section provides a comparison between our cycle and similar cycles in the literature.
\smallskip

\paragraph{Background:}
We use the term heteroclinic cycle as in \cite{AshwinChossat} where precise definitions and further detail can be found. We assume the reader is somewhat familiar with robust heteroclinic cycles in a symmetric context, for a comprehensive overview we refer to Krupa \cite{Krupa97}. In what follows we consider dynamics induced by an ODE
\begin{equation}\label{eq:equation}
\dot{x}=f(x),
\end{equation}
where $x \in \R^n$ and $f$ is smooth and $\Gamma$-equivariant for some finite group $\Gamma \subset O(n)$.

Given two hyperbolic equilibria $\xi_i$ and $\xi_j$ of system \eqref{eq:equation}, a connecting trajectory between them exists in $W^u(\xi_i) \cap W^s(\xi_j)$ if this intersection is non-empty. A heteroclinic cycle is a sequence of such connecting trajectories among a set of finitely many distinct equilibria $\xi_1, \hdots, \xi_m$ such that $\xi_{m+1}=\xi_1$. The heteroclinic cycle is the union of the equilibria and the connections.
In what follows $C_{ij}=W^u(\xi_i) \cap W^s(\xi_j)$ denotes the set of trajectories connecting two equilibria $\xi_i$ and $\xi_j$, also called a connection. 

Several notions of stability and instability are appropriate for heteroclinic cycles. We use that of completely unstable as in \cite{KrupaMelbourne95b}. Several other notions can be found in \cite{Melbourne1991}, \cite{Brannath}, \cite{KS} or \cite{PodviginaAshwin2011}. The definition of simple, pseudo-simple and quasi-simple can be found in \cite{KrupaMelbourne2004}, \cite{Podvigina2012} and \cite{PodviginaLohse} for the first, \cite{PodviginaChossat} and \cite{ChossatLohsePodvigina} for the second, and \cite{GarridoDaSilvaCastro} for the last.

\section{A hybrid heteroclinic cycle}
Our example consists in the following modification of a vector field generating a $(B_2^+,B_2^+)$ network\footnote{In \cite{CastroLohse2014} another definition of heteroclinic cycle and network is used which is why the heteroclinic object is called a network. In the context of the present article it is a heteroclinic cycle.} from \cite{CastroLohse2014},
\begin{equation}\label{eq:example}
\left\{ \begin{array}{l}
\dot{x}_1 = x_1 + \sum_{i=1}^4 b_{1i}x_i^2 + c_1x_1^3  + d_1x_3x_4\\
\dot{x}_2 = x_2 + x_2\sum_{i=1}^4 b_{2i}x_i^2 + d_2x_1x_2 \\
\dot{x}_3 = x_3 + x_3\sum_{i=1}^4 b_{3i}x_i^2 + c_3x_3^2x_4 + d_3x_1x_3 + \epsilon_3x_4\\
\dot{x}_4 = x_4 + x_4\sum_{i=1}^4 b_{4i}x_i^2 + c_4x_3x_4^2 + d_4x_1x_4 + \epsilon_4x_3
\end{array} \right. ,
\end{equation}
where all constants are real and chosen conveniently below.  This vector field is equivariant under the action of the group $\Gamma \cong \Z_2^2$ generated by
\begin{eqnarray*}
\kappa_2 . (x_1,x_2,x_3,x_4) = (x_1,-x_2,x_3,x_4) & \mbox{ and   } & \kappa_{34} . (x_1,x_2,x_3,x_4) = (x_1,x_2,-x_3,-x_4).
\end{eqnarray*}
The isotypic decomposition of $\R^4$ with respect to $\Gamma$ is
\begin{align}\label{iso-decomp}
\R^4=L_1 \oplus L_2 \oplus P_{34},
\end{align}
where $L_i$ is the $i$-th coordinate axis and $P_{ij}=L_i \oplus L_j$.

\begin{proposition}\label{th:counter-example}
System \eqref{eq:example} with coefficients fulfilling the conditions in Table~\ref{table:signs} possesses a heteroclinic cycle, $C$, between two equilibria $\xi_a$, $\xi_b \in L_1$ with a one-dimensional connection $C_{ab} \subset P_{12}$ and a two-dimensional connection $C_{ba} \subset S_{134}:=L_1 \oplus P_{34}$.
\end{proposition}

\begin{proof}
We first prove the existence of a heteroclinic cycle in the special case $d_1=\epsilon_3=\epsilon_4=0$ and then argue that it persists under generic $\Z_2^2$-equivariant perturbations. To this end, we choose coefficients in an open set of the other parameters appearing in \eqref{eq:example}  such that
\begin{enumerate}
	\item[(i)]  there are two equilibria $\xi_a=(x_a,0,0,0),\ \xi_b=(x_b,0,0,0) \in L_1$ with $x_a < 0 < x_b$;
	
	\item[(ii)] $\xi_a$ is a saddle in $P_{12}$ and a sink in $S_{134}$, and vice versa for $\xi_b$;
	
	\item[(iii)] there is a connection $C_{ab} \subset P_{12}$ from $\xi_a$ to $\xi_b$;
		
	\item[(iv)] there is a continuum of connections $C_{ba} \subset S_{134}$  from $\xi_b$ to $\xi_a$.
\end{enumerate}
We proceed step by step, collecting all conditions we obtain to satisfy (i)-(iv) in Table \ref{table:signs}.
\begin{enumerate}
	\item[(i)]  Solving for equilibria in $L_1$ yields $\dot x_1 = 0 \; \Leftrightarrow \; 1+b_{11}x_1+c_1x_1^2=0.$
		Choosing $b_{11}\neq 0$ and $c_1<0$ we obtain $\xi_a$ and $\xi_b$ with $x_a \neq -x_b$, where
		\begin{equation}\label{eq:zeros}
		x_a:= \frac{1}{2c_1} \left( {-b_{11} + \sqrt{b_{11}^2-4c_1}} \right)<0 \quad \text{and} \quad x_b:= \frac{1}{2c_1} \left( {-b_{11} - \sqrt{b_{11}^2-4c_1}} \right)>0.
		\end{equation}
	\item[(ii)]  At these equilibria the Jacobian matrix is diagonal with the following entries:
\begin{align*}
 &b_{11}x_{a/b}+2c_1x_{a/b}^2, \; 1+b_{21}x_{a/b}^2+d_2x_{a/b},\; 1+b_{31}x_{a/b}^2+d_3x_{a/b}, \; 1+b_{41}x_{a/b}^2+d_4x_{a/b}.
\end{align*}
We investigate the signs of its eigenvalues. Using the fact that the coordinates $x_a$ and $x_b$ of the equilibria satisfy \eqref{eq:zeros}, the eigenvalue along $L_1$ is, respectively, 
$$
b_{11}x_a+2c_1x_a^2=x_a\sqrt{b_{11}^2-4c_1} < 0 \;\;\; \mbox{ and  } \;\;\; b_{11}x_b+2c_1x_b^2=-x_b\sqrt{b_{11}^2-4c_1} < 0.
$$
Since $x_a < 0 < x_b$, by choosing $d_2 - \frac{b_{21}}{c_1} b_{11} < 0$ and large in absolute value, we achieve the desired signs:
$$
1-\frac{b_{21}}{c_1} + \left( d_2 - \frac{b_{21}}{c_1} b_{11} \right) x_{a} > 0  \;\;\; \mbox{ and  } \;\;\; 1-\frac{b_{21}}{c_1} + \left( d_2 - \frac{b_{21}}{c_1} b_{11} \right) x_{b} < 0.
$$

Analogous calculations for the eigenvalues along $L_3$ and $L_4$ show that if $d_3 - \frac{b_{31}}{c_1} b_{11} > 0$ and $d_4 - \frac{b_{41}}{c_1} b_{11} > 0$ and large in absolute value, both eigenvalues are negative at $\xi_a$ and positive at $\xi_b$.

	\item[(iii)] We look at the $x_1$- and $x_2$-nullclines in $P_{12}$. Because of the symmetry we restrict to $x_2>0$. We have $\dot x_1=0$ if and only if $x_2^2 = -\frac{1}{b_{12}}\left( x_1 + b_{11}x_1^2 + c_1x_1^3  \right) =:g_1(x_1)$. The zeros of $g_1$ are at $x_a$,  $x_b$ and $0$ and choosing $b_{12}>0$ we get an $x_1$-nullcline as in Figure \ref{fig:nullclines} (left). Similarly we have $\dot x_2=0$ if and only if $x_2^2 =  -\frac{1}{b_{22}}\left( 1 + b_{21}x_1^2 + d_2x_1 \right)=:g_2(x_1)$. The zeros of $g_2$ are at
	$$x_{\pm} :=  \frac{-d_2 \pm \sqrt{d_2^2 - 4b_{21}}}{2b_{21}}.$$
Fixing $b_{21}>0$ and $b_{22}<0$ we choose $d_2<0$ sufficiently large in absolute value to ensure $0<x_- < x_b < x_+$. This leads to an $x_2$-nullcline as in Figure \ref{fig:nullclines} (left).  Due to the asymptotic behaviour of $g_1$ and $g_2$ the nullclines do not intersect and there are no equilibria in $P_{12}$ off the coordinate axes. The existence of a connection $C_{ab}$ then follows from the Poincar\'{e}-Bendixson theorem. Note that the location of the nullclines in $P_{12}$ ensures that the unstable manifold $W^u(\xi_a)$ is bounded: it starts out with $\dot{x}_1=0$ and $\dot{x}_2>0$. Between the left segments of the dotted and solid curves we have $\dot{x}_1>0$ and $\dot{x}_2>0$, so eventually $W^u(\xi_a)$ must cross the $x_2$-nullcline. By the same reasoning it must cross the right segment of the $x_1$-nullcline, but cannot cross the right segment of the $x_2$-nullcline, therefore converging to $\xi_b$.

  \begin{figure}
     \centering
     \begin{subfigure}[b]{0.45\textwidth}
         \centering
\includegraphics[width=\textwidth]{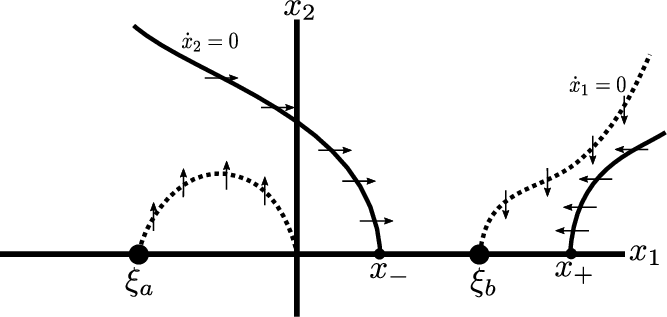}
     \end{subfigure}
  \qquad
     \begin{subfigure}[b]{0.45\textwidth}
         \centering
\includegraphics[width=\textwidth]{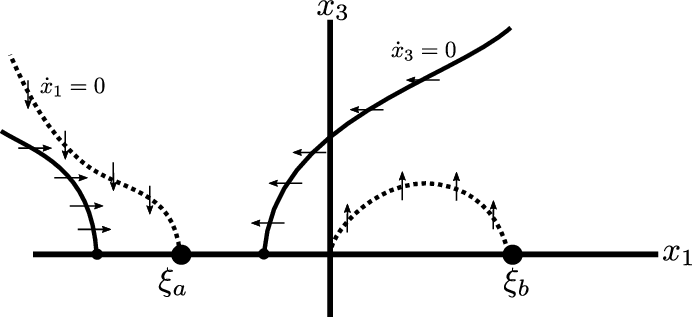}
     \end{subfigure}
      \caption{Nullclines and dynamics in $P_{12}$ (left) and $P_{13}$ (right): $\dot{x}_1=0$ dotted, $\dot{x}_2=0$ (resp.\ $\dot{x}_3=0$) solid.
\label{fig:nullclines}}
\end{figure}
 
	\item[(iv)]  The planes $P_{13}$ and $P_{14}$ are flow-invariant if $\epsilon_3=\epsilon_4=0$. We ensure that in this case there are connections from $\xi_b$ to $\xi_a$ in the same way as in (iii) by setting $b_{13}, b_{14}<0$ (for the $x_1$-nullclines) as well as $b_{33}, b_{44} <0$ and $b_{31}, b_{41}>0$ together with $d_3, d_4>0$ sufficiently large (for the $x_3$- and $x_4$-nullclines, respectively), leading to Figure \ref{fig:nullclines} (right).

Next, we observe that the set of points $Z$ in $S_{134}$ where $\dot{x}_1=0$ is determined by $x_1 + b_{11}x_1^2 + c_1x_1^3=-b_{13}x_3^2 - b_{14}x_4^2$ if $d_1=0$. In planes of constant $x_1$ this is an elipse, see Figure~\ref{fig:zeros}. Outside of $Z$ we have $\dot{x}_1<0$ (inside we have $\dot{x}_1>0$) provided $b_{13}, b_{14} < 0$, so that trajectories on $W^u(\xi_b)$, which is tangent to the affine plane spanned by $x_3$ and $x_4$, have decreasing $x_1$ component close to $\xi_b$ and therefore move in the direction of $\xi_a$. Since any equilibria in $S_{134}$ lie in $Z$, then all the trajectories in $W^u(\xi_b)$ approach $\xi_a$.
	
\begin{figure}[!htb]
 \centerline{\includegraphics[width=0.6\textwidth]{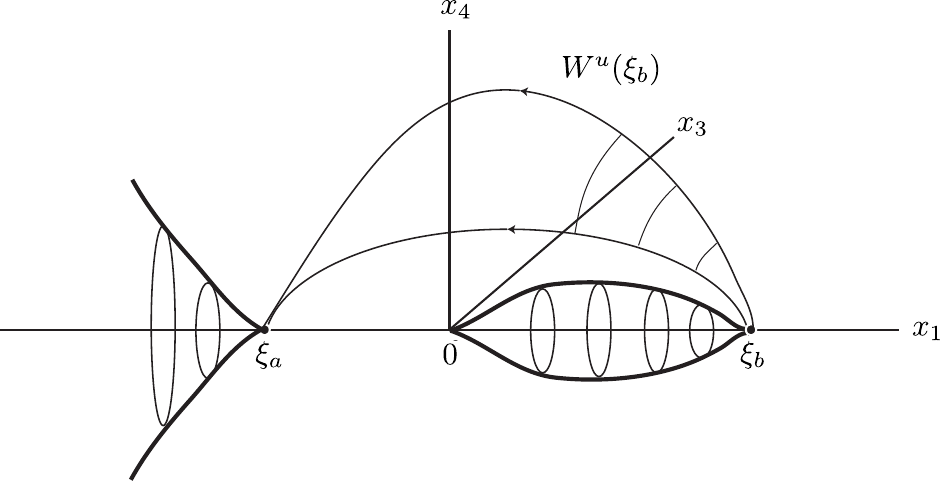}}
 \caption{The set $Z$ (thick lines) in $S_{134}$.\label{fig:zeros}}
 \end{figure}
\end{enumerate}

The conditions imposed above correspond to C1--C14 in Table~\ref{table:signs}. It is clear that these can be satisfied simultaneously and define an open set in the space of all coefficients. This finishes the proof of (i)-(iv).

 Now consider the perturbed system where $d_1, \epsilon_3, \epsilon_4 \neq 0$: since the heteroclinic connections established above are saddle-sink connections in invariant subspaces that persist under $\Gamma$-equivariant perturbations ($P_{12}$ and $S_{134}$, respectively), the cycle persists under sufficiently small such perturbations as well.  
\end{proof}

\begin{table}
\begin{center}
\begin{tabular}{| c | c | l || c | c | l |}
\hline 
 & Inequality & ensuring & & Inequality & ensuring  \\
 \hline 
C1 & $b_{11} \neq 0$ & (i) & C8 & $d_2-\frac{b_{21}}{c_1}b_{11} < 0$, large in abs. value & (ii)  \\
 C2 & $c_1< 0$ & & C9 & $d_3-\frac{b_{31}}{c_1}b_{11} > 0$, large & \\
 \cline{1-3}
 C3 & $b_{12} >0$ & (iii)&  C10 & $d_4-\frac{b_{41}}{c_1}b_{11} > 0$, large & \\
 \cline{4-6}
C4 & $b_{21} >0$ & & C11 & $b_{13}, b_{14} < 0$ & (iv) \\
C5 & $b_{22} <0$ & & C12 & $b_{31}, b_{41} > 0$ & \\
C6 & $d_2 <0$, large in abs. value & & C13 & $b_{33}, b_{44} < 0$ &   \\
C7 & $d_2^2-4b_{21} >0$ & & C14 & $d_3, d_4 > 0$, large &  \\
 \hline
\end{tabular}
\end{center}
\caption{\small{List of conditions imposed on the coefficients of \eqref{eq:example} in the construction of the vector field and the proof of Proposition~\ref{th:counter-example}. The first column assigns a label to each condition and the third column states the purpose of the restriction. The middle column contains the conditions on the coefficients.} \label{table:signs}}
\end{table}

When $\epsilon_3=\epsilon_4=d_1=0$, we define two subcycles, $C_3$ and $C_4$, of the heteroclinic cycle $C$. Both subcycles have only one-dimensional connections $[\xi_a \longrightarrow \xi_b \longrightarrow \xi_a]$. The cycle $C_i$ is such that the connection $[\xi_b \longrightarrow \xi_a]$ in $P_{1i}$ for $i=3,4$. Each subcycle $C_i$ is contained in $S_{12i}$, $i=3,4$.

\begin{lemma}\label{lem:instability}
The two subcycles $C_i$, $i=3,4$, are completely unstable in $S_{12i}$, for each $i=3,4$. 
\end{lemma}

\begin{proof}
The eigenvalues at each equilibrium in $S_{12i}$, $i=3,4$, are
\begin{itemize}
	\item  at $\xi_a$: $e_a = 1+b_{21}x_a^2+d_2 x_a>0$, in the direction of $L_2$, and $-c_{ai} = 1+b_{i1}x_a^2+d_i x_a<0$, in the direction of $L_i$;
	
	\item  at $\xi_b$: $-c_b  = 1+b_{21}x_b^2+d_2 x_b < 0$, in the direction of $L_2$, and $e_{bi} = 1+b_{i1}x_b^2+d_i x_b > 0$, in the direction of $L_i$.
\end{itemize}
We use Theorem 2.4 in \cite{KrupaMelbourne95b} to prove our claim. We check the hypotheses (S1--S3) in \cite{KrupaMelbourne95b} since (S4) is satisfied in three-dimensional space. These are, in our context:
\begin{itemize}
\item[(S1)] there exists a subgroup $\Sigma_{a} = \langle \kappa_2\rangle$, whose fixed-point space $P_{12}$ satisfies $W^u(\xi_a) \cap P_{12} \subset W^s(\xi_b)$ and $\xi_b$ is a sink in $P_{12}$; also, there exists a subgroup $\Sigma_{b} = \langle\kappa_{34}\rangle$ such that $P_{1i}=\mbox{Fix}\langle\kappa_{34}\rangle \subset S_{12i}$ where $\xi_a$ is a sink and $W^u(\xi_b) \cap P_{1i} \subset W^s(\xi_a)$, for $i=3,4$.

\item[(S2)] the eigenspaces corresponding to $-c_{ai}$ and $e_{bi}$ are in the same isotypic component since the eigenspaces coincide with $L_i$, $i=3,4$; analogously, the eigenspaces of $-c_b$ and $e_a$ coincide with $L_2$ and are therefore in the same isotypic component (see the isotypic decomposition given in \eqref{eq:isotypic}).

\item[(S3)] dim$(W^u(\xi_a) \cap P_{12}) =$ dim$(W^u(\xi_b) \cap P_{1i}) =1$ in $S_{12i}$ for $i=3,4$ (recall that dim $(N(\Sigma)/\Sigma)=0$ when $\Sigma$ is finite as is the case here).
\end{itemize}
According to Theorem 2.4 in \cite{KrupaMelbourne95b}, the subcycle $[\xi_a \longrightarrow \xi_b \longrightarrow \xi_a]$ contained in $S_{12i}$, $i=3,4$, is completely unstable if $\rho < 1$, where $\rho = \rho_{ai} \rho_{bi}$ and

$$
\rho_{ai}=\frac{c_{ai}}{e_a} \mbox{  and    } \rho_{bi}=\frac{c_b}{e_{bi}}; i=3,4.
$$
We have, in $S_{12i}$, for $i=3,4$, 
\begin{eqnarray*}
\rho < 1 & \Leftrightarrow & (1+b_{21}x_a^2+d_2 x_a)(1+b_{i1}x_b^2+d_i x_b) - (1+b_{i1}x_a^2+d_i x_a)(1+b_{21}x_b^2+d_2 x_b) > 0
\end{eqnarray*}
which, taking into account that $x_b-x_a >0$, $x_b+x_a = -b_{11}/c_1$, and $x_a x_b = 1/c_1 < 0$, becomes
\begin{eqnarray*}
 &  & b_{i1}(x_b^2-x_a^2) + d_i (x_b - x_a) -b_{21}(x_b^2-x_a^2) - b_{21}d_i(x_b-x_a)x_a x_b - d_2(x_b-x_a) + b_{i1}d_2(x_b-x_a)x_a x_b >0\\
 & \Leftrightarrow & b_{i1}(x_b+x_a)+d_i - b_{21}(x_b+x_a) - b_{21}d_ix_a x_b - d_2 + b_{i1}d_2x_a x_b > 0  \\
 & \Leftrightarrow & (d_i - b_{i1}\frac{b_{11}}{c_1}) + (b_{21}\frac{b_{11}}{c_1} -d_2) - \frac{1}{c_1} (b_{21}d_i-b_{i1}d_2) > 0.
\end{eqnarray*}
Using Table~\ref{table:signs}, the first term above is positive if $i=3$ by C9 (C10 if $i=4$), as is the second by C8. By C2, $c_1<0$ and by C6, $d_2<0$, with all other constants positive, so that the remaining term is positive. Hence, $\rho < 1$.

\end{proof}

\begin{remark}
The heteroclinic cycle $C$ constructed above has the following properties:
\begin{itemize}
 \item All trajectories in the two-dimensional connection $C_{ba}$ have the same isotropy type.
\item The heteroclinic cycle contains all unstable manifolds of its equilibria.
 \item All eigenspaces tangent to (incoming or outgoing) connections at an equilibrium are fully contained in a single $\Gamma$-isotypic component.
\end{itemize}
Note that together with (ii) from the proof above this means that $C$ satisfies hypothesis (H) in \cite{AshwinChossat}, whereas by Lemma~\ref{lem:instability} it does not satisfy hypothesis (3). Since some conditions in Table~\ref{table:signs} are merely sufficient, further analysis of this example may contribute to a better understanding of a conjecture in \cite{AshwinChossat} concerning the stability of subcycles associated with the strongest expanding eigenvalues.
\end{remark}

\section{Comparison with similar cycles}

There are several instances in the literature where similar heteroclinic structures in $\R^4$ have been studied. Some of the earliest discussions of cycles with one-dimensional connections between two equilibria can be found in \cite{Brannath} and \cite{Melbourne1991}, where conditions for essential and relative asymptotic stability are derived.

In the appendix of \cite{CastroLohse2014}  a so-called $(B_2^+, B_2^+)$ network is constructed: in a $\Z_2^3$-equivariant system it consists of two equilibria on the same coordinate axis and two simple cycles of type B connecting them in different coordinate planes. The presence of a two-dimensional connection is not discussed, even though it seems to be often present in one direction. For the one-dimensional simple cycles, a detailed stability analysis is given using the stability index defined by \cite{PodviginaAshwin2011}, allowing to draw conclusions about the stability of the full network as well. 

In \cite{CastroLohse2016} the symmetry of the system in \cite{CastroLohse2014} is partially broken so that the simple cycles making up the heteroclinic structure are of type A instead of B, making it an $(A_2^+, A_2^+)$ network in a $\Z_2^2$-equivariant system. Again, a stability index analysis is possible and the attraction configurations turn out to be not nearly as rich as those for the $(B_2^+, B_2^+)$ network.

In section 4.1 of \cite{PodviginaChossat} a pseudo-simple cycle between two equilibria in a $\D_3$-equivariant system is studied: connections are contained in planes and one of them has several non-trivial symmetric copies. In contrast to simple cycles, pseudo-simple cycles are characterized by at least one two-dimensional $\Delta_j$-isotypic component, where $\Delta_j$ is the subgroup fixing the one-dimensional space containing the equilibrium $\xi_j$.

Another example is found in \cite{ChossatLohsePodvigina}: in a $\D_n \times \Z_2$-equivariant system two equilibria are connected by a single trajectory in one direction, and by a two-dimensional intersection of invariant manifolds in the other. There is a pseudo-simple subcycle with some similarities to $C_3$ and $C_4$ in our example. The whole cycle is shown to be asymptotically stable, while the pseudo-simple subcycle attracts (at least) a set of positive measure.

Finally, a generalization of the stability index approach in \cite{PodviginaAshwin2011} is discussed in \cite{GarridoDaSilvaCastro} for quasi-simple cycles. These are characterized by one-dimensional connections in invariant spaces of equal dimension. While this true of our cycles $C_3$ and $C_4$, we explain below why the results in \cite{GarridoDaSilvaCastro} still do not encompass our example.
\bigbreak
We now discuss how our example fits into the landscape sketched above, focusing on the unperturbed case, i.e.\ when $\epsilon_3=\epsilon_4=d_1=0$. First note that our cycle $C$ is not simple due to its two-dimensional connection $C_{ba}$. Even though the subcycles $C_3$ and $C_4$ have only one-dimensional connections, they are also not simple because of the two-dimensional $\Gamma$-isotypic component in (\ref{iso-decomp}). While such a two-dimensional component is a key feature of a pseudo-simple cycle, none of our cycles is pseudo-simple because the planes $P_{13}$ and $P_{14}$ containing the connections from $\xi_b$ to $\xi_a$ are not fixed-point spaces -- even though they are dynamically invariant. When the system is perturbed the invariance of those planes is broken, moving the cycles further away from the simple/pseudo-simple setting.

Let us take a closer look at the transitions between neighbourhoods of $\xi_a$ and $\xi_b$ for trajectories near $C$ in order to pinpoint similarities and differences to simple, pseudo-simple and quasi-simple cycles. Combining information on these global transitions with the linearization of the vector field near the equilibria to obtain a full return map is the standard method for analyzing the attraction properties of heteroclinic cycles, see \cite{KrupaMelbourne2004} for a detailed discussion. This approach also lies at the heart of the stability index results in  \cite{PodviginaAshwin2011} and \cite{GarridoDaSilvaCastro}. It is well-known that in $\R^4$ it typically suffices to consider two coordinates for the return maps, yielding the standard form found in equation (4.1) of \cite{KrupaMelbourne2004}. Additionally, the symmetry of a system imposes restrictions on its global transitions. Recall that in our example we have $C_{ab} \subset P_{12}=\Fix(\Sigma_a)$ with $\Sigma_a:=\{\mathbbm{1}, \kappa_{34} \}$, where $\xi_b$ is a sink; and $C_{ba} \subset S_{134}=\Fix(\Sigma_b)$ with $\Sigma_b:=\{\mathbbm{1}, \kappa_2 \}$, where $\xi_a$ is a sink. The respective $\Sigma_a$ and $\Sigma_b$-isotypic decompositions are
\begin{equation}\label{eq:isotypic}
\R^4 = \Fix(\Sigma_a) \oplus P_{34} = P_{12} \oplus P_{34} \quad \text{and} \quad \R^4 = \Fix(\Sigma_b) \oplus L_2 = S_{134} \oplus L_2.
\end{equation}
For the transitions from $\xi_b$ to $\xi_a$ in $C_3$ and $C_4$ the second decomposition restricts the corresponding map just as in Proposition 4.1 (ii) of \cite{KrupaMelbourne2004}, which means this connection can be viewed as simple and of type B in the sense of \cite{KrupaMelbourne2004}. In particular, a quasi-simple cycle where all transitions are like this would be accessible via the stability results in \cite{GarridoDaSilvaCastro}.

However, for the other transition from $\xi_a$ to $\xi_b$ the first decomposition above does not yield any further restrictions on the corresponding map, just as in Proposition 4.1 (i) in \cite{KrupaMelbourne2004}. Therefore, this connection may be considered as type A. In particular, it does not allow us to reduce the study of the return map to the level of transition matrices. This means that $C_3$ and $C_4$ are quasi-simple, but the results in \cite{GarridoDaSilvaCastro} do not apply to them. To our knowledge our system is the first example of a heteroclinic object with transitions of both type A and B. While in principle the return map approach could still be applied, this means that the standard reduction methods for analyzing stability mentioned above fail. 

In Theorem 1 of \cite{PodviginaChossat} it is shown that pseudo-simple cycles with symmetry group in $\SO(4)$ are generically completely unstable. Lemma \ref{lem:instability} shows that -- under the conditions in Table \ref{table:signs} -- the same is true for $C_3$ and $C_4$ in our example, even though $\Gamma \not\subset \SO(4)$. The reason for this instability, however, is clearly different from that for pseudo-simple cycles, since the proof of Theorem 1 in \cite{PodviginaChossat} crucially depends on the fact that a pseudo-simple cycle always has at least two equilibria where two eigenvalues of the linearization are equal. This is a consequence of the symmetry of pseudo-simple cycles, and not the case in our example. Numerical experiments seem to indicate that in our system even $C$ is generically completely unstable.
\bigbreak
In conclusion, we have described an elementary example of a heteroclinic cycle between two equilibria in $\R^4$ that seems to fit in none of the established categories of low-dimensional heteroclinic objects. Its attraction and stability is not accessible using standard techniques for return map reduction or stability indices. Our comparison with similar heteroclinic objects in the literature further illustrates the subtle influence of symmetry on equivariant dynamics.
\bigbreak

\paragraph{Acknowledgements:}
The authors are grateful to P.\ Ashwin for helpful comments.

The first author was partially supported by Centro de Matemática da Universidade do Porto
(CMUP), financed by national funds through FCT - Fundação para a Ciência
e a Tecnologia, I.P., under the project UIDB/00144/2020.

Both authors benefitted from DAAD-CRUP funding through ``A\c{c}\~ao Integrada Luso-Alem\~a A10/17'', which on the German side is funded as project 57338573 PPP Portugal 2017 by the German Academic Exchange Service (DAAD), sponsored by the Federal Ministry of Education and Research (BMBF).

Declarations of interest: none.

\end{document}